\documentclass[a4paper]{llncs}
\usepackage[centertags]{amsmath} 
\usepackage{amssymb}
\usepackage{graphicx}
\usepackage{bbm}

\usepackage{ifthen}
\newcommand{\isdraft}{\boolean{true}} 
\renewcommand{\isdraft}{\boolean{false}} 
\ifthenelse{\isdraft}{
    \usepackage[color]{showkeys}
    \usepackage{xcolor}
}{}
\newcommand{\markupdraft}[2]{
    \ifthenelse{\equal{#1}{display}}{#2}{}
    \ifthenelse{\equal{#1}{color}}{\color{#2}}{}
}

\newcommand{\newcolored}[3][]{{\markupdraft{color}{#2}#3}
    \ifthenelse{\equal{#1}{}}{}{\markupdraft{display}{{\color{yellow!70!black}[#1]}}}} 

\providecommand{\del}[2][]{{\markupdraft{display}{{\color{red!20!yellow}[rmed: "#2"[#1]]}}}} 

\providecommand{\todo}[2][]{\markupdraft{display}{{\color{red}\noindent++TODO: #2 {\color{yellow}(#1)}++}}}
\ifthenelse{\isdraft}{}{\renewcommand{\markupdraft}[2]{}}

\newcommand{\mc}[2]{\newcommand{#1}{\ensuremath{#2}}}
\newcommand{\bs}[1]{\boldsymbol{#1}}
\mc{\R}{\mathbb{R}}
\mc{\N}{\mathbb{N}}
\newcommand{\ei}[1][k]{\ensuremath{\bs{e}_{#1}}}
\mc{\x}{\bs{x}}
\mc{\X}{\bs{X}}
\newcommand{\Xt}[1][t]{\ensuremath{\X_{#1}}}
\newcommand{\Xtt}{\Xt[t+1]}
\mc{\Xfea}{\mathcal{X}_{\textrm{feasible}}}

\newcommand{\st}[1][t]{\ensuremath{\sigma_{#1}}}
\newcommand{\stt}{\st[t+1]}
\mc{\y}{\bs{y}}
\mc{\Y}{\bs{Y}}
\newcommand{\Yti}[1][i]{\ensuremath{\Y_t^{#1}}}
\mc{\M}{\bs{M}}
\newcommand{\Mtij}[1][j]{\M_t^{i,#1}}
\newcommand{\Mti}[1][i]{\M_t^{#1}}
\newcommand{\Mtstar}[1][t]{\ensuremath{\M_{#1}^\star}}
\mc{\dist}{H}
\newcommand{\distilde}[1][\delta]{\tilde{\dist}_{#1}}
\newcommand{\distar}[1][\delta]{\distilde[#1]^{\star}}
\newcommand{\den}{h}
\newcommand{\dentilde}[1][\delta]{\tilde{\den}_{#1}}
\newcommand{\denstar}[1][\delta]{\dentilde[#1]^{\star}}
\newcommand{\dt}[1][t]{\delta_{#1}}
\newcommand{\dtt}{\dt[t+1]}
\newcommand{\dmarkov}{(\dt)_{t\in \N}}
\newcommand{\seqmtij}{(\Mtij)_{i \in [1..\lambda], (j,t)\in \N^2}}
\newcommand{\W}{\mathbf{W}}
\newcommand{\Wt}[1][t]{\W_{#1}}
\newcommand{\E}{\mathbf{E}}
\newcommand{\G}{\mathcal{G}}
\newcommand{\Gstar}{\G^\star}
\newcommand{\Q}{\mathbf{Q}}
\newcommand{\C}{\mathbf{C}}
\newcommand{\cdfk}[1][k]{F_{#1,\delta}}

\newcommand{\cdfkinv}[1][k]{F_{#1,\delta}^{-1}}

\newcommand{\ddomain}{\R_+}
\newcommand{\eqd}{\overset{(d)}{=}}
\newcommand{\Id}[1][n]{\mathbf{I}_{#1}}
\newcommand{\Nl}{\mathcal{N}}
\newcommand{\Nlmult}{\Nl(\bs{0},\Id)}
\newcommand{\Nlambda}[1][\lambda]{\Nl_{#1:\lambda}}
\newcommand{\Ld}[1][\delta]{L_{#1}}
\newcommand{\ind}{\mathbbm{1}}
\newcommand{\borel}{\mathcal{B}}
\newcommand{\U}{\bs{U}}
\newcommand{\Ldstar}[1][v]{L_{\delta, #1}^\star}
\newcommand{\Va}[1][\alpha]{V_{#1}}

\newcommand{\mplus}{\mu_{+}}

\newcommand{\B}{\bs{B}}
\newcommand{\Ntij}{\bs{N}_t^{i,j}}
\newcommand{\Nti}{\bs{N}_t^i}
\newcommand{\Ntstar}{\bs{N}_t^\star}
\newcommand{\ekb}[1][k]{\bar{\bs{e}}_{#1}}

\newcommand{\ekbp}[1][k]{\bar{\bs{e}}_{#1}'}

\newcommand{\ekp}[1][k]{\bs{e}_{#1}'}

\newcommand{\normb}[1]{\| #1 \|_{-}}

\newcommand{\denstartwo}{\den_{\delta,2}^\star}
\renewcommand{\u}{\bs{u}}

\DeclareMathOperator{\argmax}{argmax}
\DeclareMathOperator{\Var}{Var}

\newenvironment{myproof}{\begin{proof}}{\qed \end{proof}} 

\newcounter{enumroman}
\newenvironment{romanitems}{\begin{list}{\bfseries(\roman{enumroman})\hfill}{\usecounter{enumroman}
\setlength{\labelwidth}{\leftmargin}\addtolength{\labelwidth}{-1\labelsep}
\topsep=0mm plus 2pt\itemsep=0mm plus 1pt\parsep=0mm\itemsep=0mm plus
1pt\itemindent=0mm}}{\end{list}}

\title{           A Generalized Markov-Chain Modelling Approach to $(1,\lambda )$-ES Linear Optimization: Technical Report}
\author{          Alexandre Chotard\,\inst{1} \and Martin Hole\v{n}a\,\inst{2}} 
\institute{       
                   INRIA Saclay-Ile-de-France, LRI, \email{alexandre.chotard@lri.fr}
                  University Paris-Sud, France \and Institute of Computer Science, Academy of Sciences, Pod vod\'arenskou v\v{e}\v{z}\'{\i} 2, Prague, Czech Republic,
                  \email{martin@cs.cas.cz}}

\begin{document}

\maketitle

\begin{abstract}
Several recent publications investigated Markov-chain modelling of linear optimization by a $(1,\lambda )$-ES, considering both unconstrained and linearly constrained optimization, and both constant and varying step size. All of them assume normality of the involved random steps, and while this is consistent with a black-box scenario, information on the function to be optimized (e.g. separability) may be exploited by the use of another distribution. 
The objective of our contribution is to complement previous studies realized with normal steps, and to give sufficient conditions on the distribution of the random steps for the success of a constant step-size $(1,\lambda)$-ES on the simple problem of a linear function with a linear constraint. The decomposition of a multidimensional distribution into its marginals and the copula combining them is applied to the new distributional assumptions, particular attention being paid to distributions with Archimedean copulas.
\end{abstract}    

\section{Introduction} 
Evolution Strategies (ES) are Derivative Free Optimization (DFO) methods, and as such are suited for the optimization of numerical problems in a black-box context, where the algorithm has no information on the function $f$ it optimizes (e.g. existence of gradient) and can only query the function's values. In such a context, it is natural to assume normality of the random steps, as the normal distribution has maximum entropy for given mean and variance, meaning that it is the most general assumption one can make without the use of additional information on $f$. However such additional information may be available, and then using normal steps may not be optimal. Cases where different distributions have been studied include so-called Fast Evolution Strategies \cite{yao1997fast} or SNES \cite{bbob-snes,Schaul2011snes} which exploits the separability of $f$, or heavy-tail distributions on multimodal problems \cite{hansen2006wdh,Schaul2011snes}.

In several recent publications \cite{arnold2011behaviour,arnold2012behaviour,cahPPSN12,cah2014linearconstraint}, attention has been paid to Markov-chain modelling of linear optimization by a $(1,\lambda )$-ES, i.e. by an evolution strategy in which $\lambda $ children are generated from a single parent  $\X\in\mathbb{R}^n$ by adding normally distributed $n$-dimensional random steps $\M$, \small 
\begin{gather}
\label{eq:esintro}
                  \X \leftarrow \X+ \sigma \C^{\frac{1}{2}} \M, \text{ where }\M\sim \Nlmult.
\end{gather}      
\normalsize Here, $\sigma $ is called step size, $\C$ is a covariance matrix, and $\Nlmult$ denotes the $n$-dimensional standard normal distribution with zero mean and covariance matrix identity. The best among the $\lambda $ children, i.e. the one with the highest fitness, becomes the parent of the next generation, and the step-size $\sigma$ and the covariance matrix $\C$ may then be adapted to increase the probability of sampling better children. In this paper we relax the normality assumption of the movement $\M$ to a more general distribution $\dist$.

The linear function models a situation where the step-size is relatively small compared to the distance towards a local optimum. This is a simple problem that must be solved by any effective evolution strategy by diverging with positive increments of $\nabla f . \M$. \del{As ES are comparison based algorithm they are invariant to the composition of the fitness function $f$ by any strictly increasing function, $f$ could be composed by a function without derivative and with many discontinuities without any impact on our study.} This unconstrained case was studied in \cite{cahPPSN12} for normal steps with cumulative step-size adaptation (the step-size adaptation mechanism in CMA-ES \cite{cmaes}).

Linear constraints naturally arise in real-world problems (e.g. need for positive values, box constraints) and also model a step-size relatively small compared to the curvature of the constraint.   Many techniques to handle constraints in randomised algorithms have been proposed (see \cite{CoelloCoello:2008}). In this paper we focus on the resampling method, which consists in resampling any unfeasible candidate until a feasible one is sampled. We chose this method as it makes the algorithm easier to study, and is consistent with the previous studies assuming normal steps \cite{arnold2008behaviour,arnold2011behaviour,arnold2012behaviour,cah2014linearconstraint}, studying constant step-size, self adaptation and cumulative step-size adaptation mechanisms (with fixed covariance matrix).

Our aim is to study the $(1,\lambda)$-ES with constant step-size, constant covariance matrix and random steps with a general absolutely continuous distribution $\dist$ optimizing a linear function under a linear constraint handled through resampling. We want to extend the results obtained in \cite{arnold2011behaviour,cah2014linearconstraint} using the theory of Markov chains. It is our hope that such results will help in designing new algorithms using information on the objective function to make non-normal steps. We pay a special attention to distributions with Archimedean copulas, which are a particularly well transparent alternative to the normal distribution. Such distributions have been recently considered in the Estimation of Distribution Algorithms \cite{cuestainfante10bivariate,wang12copula}, continuing the trend of using copulas in that kind of evolutionary optimization algorithms                  \cite{salinasgutierez09using}.

In the next section, the basic setting for modelling the considered evolutionary optimization task is formally defined. In Section~\ref{sc:dist}, the distributions of the feasible steps and of the selected steps are linked to the distribution of the random steps, and another way to sample them is provided. In Section~\ref{sec:markov}, it is shown that, under some conditions on the distribution of the random steps, the normalized distance to the constraint is a ergodic Markov chain, and a law of large numbers for Markov chains is applied. Finally, Section~\ref{sec:examples} gives properties on the distribution of the random steps under which some of the aforementioned conditions are verified.

\subsection*{Notations}
For $(a,b) \in \N^2$ with $a < b$, $[a..b]$ denotes the set of integers $i$ such that $a \leq i \leq b$. For $X$ and $Y$ two random vectors, $X \eqd Y$ denotes that these variables are equal in distribution, $X \overset{a.s.}{\rightarrow} Y$ and $X \overset{\mathcal{P}}{\rightarrow} Y$ denote, respectively, almost sure convergence and convergence in probability. For $(\x,\y) \in \R^n$, $\x.\y$ denotes the scalar product between the vectors $\x$ and $\y$, and for $i \in [1..n]$, $[\x]_i$ denotes the $i^{\textrm{th}}$ coordinate of $\x$. For $A$ a subset of $\R^n$,  $\ind_A$ denotes the indicator function of $A$. For $\mathcal{X}$ a topological set, $\borel(\mathcal{X})$ denotes the Borel algebra on $\mathcal{X}$.

\section{Problem setting and algorithm definition}

Throughout this paper, we study a $(1,\lambda)$-ES optimizing a linear function $f : \R^n \rightarrow \R$ where $\lambda \geq 2$ and $n \geq 2$, with a linear constraint $g : \R^n \rightarrow \R$, handling the constraint by resampling unfeasible solutions until a feasible solution is sampled.

Take $(\ei[k])_{k \in [1..n]}$ a orthonormal basis of $\R^n$. We may assume $\nabla f$ to be normalized as the behaviour of an ES is invariant to the composition of the objective function by a strictly increasing function (e.g. $h : x \mapsto x/\|\nabla f\|$), and the same holds for $\nabla g$ since our constraint handling method depends only on the inequality $g(\x) \leq 0$ which is invariant to the composition of $g$ by a homothetic transformation. Hence w.l.o.g. we assume that $\nabla f = \ei[1]$ and $\nabla g = \cos\theta \ei[1] + \sin \theta \ei[2]$ with the set of feasible solutions $\Xfea := \lbrace \x \in \R^n | g(\x) \leq 0 \rbrace$. We restrict our study to $\theta \in (0, \pi/2)$. Overall the problem reads
\begin{equation}\label{eq:pbdef}
\begin{split}
{\rm maximize}\,\,\,\, f(\x) = [\x]_{1}\,\,\,\, {\rm subject~to~}  \\ g(\x) = [\x]_1 \cos \theta + [\x]_2 \sin \theta \leq 0 \enspace.
\end{split} 
\end{equation}

\begin{figure}
\begin{center}
\includegraphics[width=0.6\textwidth]{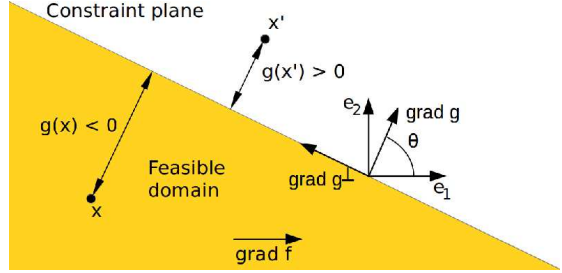}
\end{center}
\caption{Linear function with a linear constraint, in the plane spanned by $\nabla f$ and $\nabla g$, with the angle from $\nabla f$ to  $\nabla g$ equal to $\theta \in (0, \pi/2)$. The point $\x$ is at distance $g(\x)$ from the constraint hyperplan $g(\x) = 0$.}
\label{fg:problem}
\end{figure}

At iteration $t \in \N$, from a so-called parent point $\Xt \in \Xfea$ and with step-size $\st \in \R_+^*$ we sample new candidate solutions by adding to $\Xt$ a random vector $\st\Mtij$ where $\Mtij$ is called a random step and $(\Mtij)_{i \in [1..\lambda], j\in \N, t \in \N}$ is a i.i.d. sequence of random vectors with distribution $\dist$. The $i$ index stands for the $\lambda$ new samples to be generated, and the $j$ index stands for the unbounded number of samples used by the resampling. We denote $\Mti$ a feasible step, that is the first element of $(\Mtij)_{j \in \N}$ such that $\Xt + \st \Mti \in \Xfea$ (random steps are sampled until a suitable candidate is found). The $i^{\textrm{th}}$ feasible solution $\Yti$ is then 
\begin{equation}
\Yti := \Xt + \st \Mti \enspace .
\end{equation}
Then we denote $\star := \argmax_{i \in [1..\lambda]} f(\Yti)$ the index of the feasible solution maximizing the function $f$, and update the parent point
\begin{equation} \label{eq:xupdate}
\Xtt := \Yti[\star] = \Xt + \st \Mtstar \enspace ,
\end{equation}
where $\Mtstar$ is called the selected step. Then the step-size $\st$, the distribution of the random steps $\dist$ or other internal parameters may be adapted.

Following \cite{arnold2011behaviour,arnold2012behaviour,arnold2008behaviour,cah2014linearconstraint} we define $\dt$ as
\begin{equation} \label{eq:delta}
\dt := -\frac{g(\Xt)}{\st} \enspace .
\end{equation}

\section{Distribution of the feasible and selected steps}
\label{sc:dist}

In this section we link the distributions of the random vectors $\Mti$ and $\Mtstar$ to the distribution of the random steps $\Mtij$, and give another way to sample $\Mti$ and $\Mtstar$ not requiring an unbounded number of samples.
\todo{add conditional probability regarding $\delta$ in proof}

\begin{lemma} \label{lm:diststep}
Let a $(1,\lambda)$-ES optimize the problem defined in \eqref{eq:pbdef} handling constraint through resampling. Take $\dist$ the distribution of the random step $\Mtij$, and for $\delta \in \R_+^*$ denote $\Ld := \lbrace \x \in \R^n | g(\x) \leq \delta \rbrace$. Providing that $\dist$ is absolutely continuous and that $\dist(\Ld) > 0$ for all $\delta \in \R_+$, the distribution $\distilde$ of the feasible step and $\distar$ the distribution of the selected step when $\dt = \delta$ are absolutely continuous, and denoting $\den$, $\dentilde$ and $\denstar$ the probability density functions of, respectively, the random step, the feasible step $\Mti$ and the selected step $\Mtstar$ when $\dt = \delta$
\begin{equation} \label{eq:denfea}
\dentilde(\x) = \frac{\den(\x)\ind_{\Ld}(\x)}{\dist(\Ld)} \enspace ,
\end{equation}
and
\begin{align} \label{eq:denselect}
\denstar(\x) &= \lambda \dentilde(\x) \distilde((-\infty,[\x]_1) \times \R^{n-1})^{\lambda-1} \nonumber \\ 
			&=  \lambda \frac{\den(\x) \ind_{\Ld}(\x) \dist((-\infty,[\x]_1) \times \R^{n-1} \cap \Ld)^{\lambda-1}}{\dist(\Ld)^\lambda} \enspace .
\end{align}
\end{lemma}

\begin{myproof} 
Let  $\delta>0,A\in\borel(\R^n)$. Then for $t\in\mathbb{N},i=1\dots\lambda$, using the the fact that $(\Mtij)_{j \in \N}$ is a i.i.d. sequence
\begin{align*}
\distilde(A) &= \Pr(\Mti \in A | \dt = \delta) \\
&= \sum_{j\in\N} \Pr(\Mtij \in A \cap \Ld \textrm{ and } \forall k < j, \Mtij[k] \in \Ld^c | \dt = \delta) \\
&= \sum_{j \in \N} \Pr(\Mtij \in A \cap \Ld | \dt = \delta) \Pr( \forall k < j, \Mtij[k] \in \Ld^c | \dt = \delta) \\
&= \sum_{j\in\N} \dist(A \cap \Ld) (1 - \dist(\Ld))^{j} \\
&= \frac{\dist(A \cap \Ld)}{\dist(\Ld)} = \int_A \frac{h(\x) \ind_{\Ld}(\x) d\x }{\dist(\Ld)} \enspace ,
\end{align*}
which yield Eq.~\eqref{eq:denfea} and that $\distilde$ admits a density $\dentilde$ and is therefore absolutely continuous.

Since $((\Mtij)_{j \in \N})_{i\in [1..\lambda]}$ is i.i.d., $(\Mti)_{i \in [1..\lambda]}$ is i.i.d. and 
\begin{align*}
\distar(A) &= \Pr(\Mtstar \in A | \dt = \delta) \\
&= \sum_{i=1}^\lambda \Pr(\Mti \in A \textrm{ and } \forall j \in [1..\lambda] \backslash \lbrace i \rbrace,  [\Mti]_1 > [\Mti[j]]_1 | \dt = \delta) \\
&= \lambda \Pr(\Mti[1] \in A \textrm{ and } \forall j \in [2..\lambda],  [\Mti[1]]_1 > [\Mti[j]]_1 | \dt = \delta) \\
&= \lambda \int_A \dentilde(\x) \Pr( \forall j \in [2..\lambda], [\Mti[j]]_1 < [\x]_1 | \dt = \delta) d\x \\
&= \int_A \lambda  \dentilde(\x) \distilde((-\infty,[\x]_1) \times \R^{n-1})^{\lambda-1} d\x \enspace ,
\end{align*}
which shows that $\distar$ possess a density, and with \eqref{eq:denfea} yield Eq.~\eqref{eq:denselect}.
\end{myproof}

The vectors $(\Mti)_{i \in [1..\lambda]}$ and $\Mtstar$ are functions of the vectors $(\Mtij)_{i \in [1..\lambda], j\in \N}$ and of $\dt$. In the following Lemma an equivalent way to sample $\Mti$ and $\Mtstar$ is given which uses a finite number of samples. This method is useful if one wants to avoid dealing with the infinite dimension space implied by the sequence $(\Mtij)_{i \in [1..\lambda, j \in \N}$.

\begin{lemma} \label{lm:g}
Let a $(1,\lambda)$-ES optimize problem \eqref{eq:pbdef}, handling the constraint through resampling, and take $\dt$ as defined in \eqref{eq:delta}. Let $\dist$ denote the distribution of $\Mtij$ that we assume absolutely continuous, $\nabla g^\perp := -\sin\theta \ei[1] +\cos\theta \ei[2]$, $\Q$ the rotation matrix of angle $\theta$ changing $(\ei[1], \ei[2], \ldots, \ei[n])$ into $(\nabla g, \nabla g^\perp, \ldots, \ei[n])$. Take $\cdfk[1](x) := \Pr(\Mti.\nabla g \leq x | \dt = \delta)$, $\cdfk[2](x) := \Pr(\Mti.\nabla g^\perp \leq x | \dt = \delta)$ and  $\cdfk(x) := \Pr([\Mti]_k \leq x | \dt = \delta)$ for $k \in [3..n]$, the marginal cumulative distribution functions when $\dt = \delta$, and $C_\delta$ the copula of $(\Mti.\nabla g, \Mti.\nabla g^\perp, \ldots, \Mti.\ei[n])$. We define
\begin{equation} \label{eq:g}
\G : (\delta, (u_i)_{i \in [1..n]}) \in  \ddomain \times [0,1]^n \mapsto \Q \left(\begin{array}{c} \cdfkinv[1](u_1) \\ \vdots \\ \cdfkinv[n](u_n) \end{array} \right) \enspace ,
\end{equation}
\begin{equation} \label{eq:gstar}
\Gstar : (\delta, (\bs{v}_i)_{i\in [1..\lambda]}) \in \ddomain \times [0,1]^{n\lambda} \mapsto  \underset{\bs{G} \in \lbrace \G(\delta, \bs{v}_i) | i \in [1..\lambda] \rbrace}{\argmax} f(\bs{G}) \enspace .
\end{equation}
Then, if the copula $C_\delta$ is constant in regard to $\delta$, for $\Wt = (\bs{V}_{i,t})_{i\in [1..\lambda]}$ a i.i.d. sequence with $\bs{V}_{i,t} \sim C_\delta$
\begin{equation} \label{eq:gmti}
\G(\dt, \bs{V}_{i,t}) \eqd \Mti \enspace ,
\end{equation}
\begin{equation} \label{eq:gmtstar}
\Gstar(\dt, \Wt) \eqd \Mtstar \enspace .
\end{equation}
\end{lemma}

\begin{myproof}
Since $\bs{V}_{i,t} \sim \C_\delta$ 
\begin{equation*}
(\Mti.\nabla g, \Mti.\nabla g^\perp, \ldots, \Mti.\ei[n]) \eqd (\cdfkinv[1](\bs{V}_{1,t}), \cdfkinv[2](\bs{V}_{2,t}), \ldots, \cdfkinv[n](\bs{V}_{n,t})) \enspace ,
\end{equation*}
and if the function $\delta \in \ddomain \mapsto C_\delta$ is constant, then the sequence of random vectors $(\bs{V}_{i,t})_{i \in [1..\lambda], t\in \N}$ is i.i.d.. Finally by definition $\Q^{-1} \Mti = (\Mti.\nabla g, \Mti.\nabla g^\perp, \ldots, \Mti.\ei[n])$, which shows Eq.~\eqref{eq:gmti}.
Eq.~\eqref{eq:gmtstar} is a direct consequence of Eq.~\eqref{eq:gmti} and the fact that $\Mtstar = \underset{\bs{G} \in \lbrace \G(\delta, \bs{v}_i) | i \in [1..\lambda] \rbrace}{\argmax} f(\bs{G})$ (which holds as $f$ is linear).
\end{myproof}

\del{The proof can be found in \cite{chotard14TRgeneraldistributions}. The idea consists in the fact that for $\X$ a random vector with marginal cumulative distributions $(F_i)_{i \in [1..n]}$ and copula $C$, for $\U \sim C$ then $\X \eqd (F_i^{-1}([\U]_i))_{i \in [1..n]}$.
}

We may now use these results to show the divergence of the algorithm when the step-size is constant, using the theory of Markov chains \cite{markovtheory}.

\section{Divergence of the $(1,\lambda)$-ES with constant step-size}
\label{sec:markov}

Following the first part of \cite{cah2014linearconstraint}, we restrict our attention to the constant step size in the remainder of the paper, that is for all $t\in \N$ we take $\st = \sigma \in \R_+^*$.
                  

From Eq.~\eqref{eq:xupdate}, by recurrence and dividing by $t$, we see that
\begin{equation} \label{eq:xrec}
\frac{[\Xt - \Xt[0]]_1}{t} = \frac{\sigma}{t} \sum_{i=0}^{t-1} \Mtstar[i] \enspace .
\end{equation}
The latter term suggests the use of a Law of Large Numbers to show the convergence of the LHS (Left Hand Side) to a constant that we call the divergence rate. The random vectors $(\Mtstar)_{t \in \N}$ are not i.i.d. so in order to apply a Law of Large Numbers on the RHS (Right Hand Side) of the previous equation we use Markov chain theory, more precisely the fact that $(\Mtstar)_{t\in\N}$ is a function of a $(\dt, (\Mtij)_{i \in [1..\lambda], j\in \N})_{t\in\N}$ which is a geometrically ergodic Markov chain. As $(\Mtij)_{i \in [1..\lambda], j\in \N, t \in \N}$ is a i.i.d. sequence, it is a Markov chain, and the sequence $(\dt)_{t \in \N}$ is also a Markov chain as stated in the following proposition.

\begin{proposition} 
Let a $(1,\lambda)$-ES with constant step-size optimize problem \eqref{eq:pbdef}, handling the constraint through resampling, and take $\dt$ as defined in \eqref{eq:delta}. Then no matter what distribution the i.i.d. sequence $\seqmtij$ have, $\dmarkov$ is a homogeneous Markov chain and
\begin{equation} \label{eq:deltaupdate}
\dtt = \dt - g(\Mtstar) = \dt -\cos\theta[\Mtstar]_1 - \sin\theta [\Mtstar]_2 \enspace .
\end{equation}
\end{proposition}              

\begin{myproof}
By definition in \eqref{eq:delta} and since for all $t$, $\st = \sigma$,
\begin{align*}
\dtt &= -\frac{g(\Xtt)}{\stt} \\
 &= -\frac{g(\Xt) + \sigma g(\Mtstar) }{\sigma} \\
 &= \dt - g(\Mtstar) \enspace ,
\end{align*}
and as shown in \eqref{eq:denselect} the density of $\Mtstar$ is determined by $\dt$. So the distribution of $\dtt$ is determined by $\dt$, hence $\dmarkov$ is a time-homogeneous Markov chain.
\end{myproof}

We now show ergodicity of the Markov chain $\dmarkov$, which implies that the $t$-steps transition kernel (the function $A \mapsto \Pr(\dt \in A | \dt[0] = \delta)$ for $A \in \borel(\ddomain)$) converges towards a stationary measure $\pi$, generalizing Propositions 3 and 4 of \cite{cah2014linearconstraint}. 

\begin{proposition}
\label{pr:markovproperties}        
Let a $(1,\lambda)$-ES with constant step-size optimize problem \eqref{eq:pbdef}, handling the constraint through resampling. We assume that the distribution of $\Mtij$ is absolutely continuous with probability density function $\den$, and that $\den$ is continuous and strictly positive on $\R^n$. Denote $\mu_+$ the Lebesgue measure on $(\mathbb{R}_+,\mathcal{B}(\mathbb{R}_+))$, and for $\alpha>0$ take the functions $V : \delta \mapsto \delta$, $\Va : \delta \mapsto \exp(\alpha \delta)$ and $r_1 : \delta \mapsto 1$. Then $\dmarkov$ is $\mu_+$-irreducible, aperiodic and compact sets are small sets for the Markov chain.

If the following two additional conditions are fulfilled
\begin{equation}\label{eq:asint}
\E(|g(\Mtij)| \enspace | \enspace \dt = \delta) < \infty \textrm{ for all }\delta\in \ddomain \enspace, \textrm{ and} 
\end{equation}
\begin{equation} \label{eq:asineq}
\lim_{\delta \rightarrow + \infty} \E(g(\Mtstar ) | \dt = \delta) \in \R_+^* \enspace ,
\end{equation}
then $\dmarkov$ is $r_1$-ergodic and positive Harris recurrent with some invariant measure $\pi$.

Furthermore, if
\begin{equation}\label{eq:asexpint}
\E(\exp(g(\Mtij)) | \dt = \delta) < \infty \textrm{ for all } \delta \in \ddomain \enspace,
\end{equation}
then for $\alpha>0$ small enough, $\dmarkov$ is also $\Va-$geometrically ergodic.
\end{proposition}

\begin{myproof}
The probability transition kernel of $\dmarkov$ writes
\begin{align*}
 P (\delta, A) &= \int_{\R^n} \ind_{A}(\delta-g(\x)) \denstar(\x) d\x \\
 &= \int_{\R^n} \ind_{A}(\delta-g(\x)) \lambda \frac{\den(\x) \ind_{\Ld}(\x) \dist((-\infty,[\x]_1) \times \R^{n-1} \cap \Ld)^{\lambda-1}}{\dist(\Ld)^\lambda} \\
 &= \frac{\lambda}{\dist(\Ld)^\lambda}\int_{g^{-1}(A)} \den \left( \begin{array}{c} \delta - [\bs{u}]_1 \\ -[\bs{u}]_2 \\ \vdots \\ -[\bs{u}]_n \end{array} \right) \dist( (-\infty, \delta - [\bs{u}]_1) \times \R^{n-1} \cap \Ld )^{\lambda-1} d\bs{u} \enspace ,
\end{align*}
with the substitution of variables $[\bs{u}]_1 = \delta - [\x]_1$ and $[\bs{u}]_i = -[\x]_i$ for $i \in [2.. n]$. Denote $\Ldstar := (-\infty,v)\times \R^{n-1} \cap \Ld$ and $t_\delta : \bs{u} \mapsto (\delta - [\bs{u}]_1, -[\bs{u}]_2, \ldots, -[\bs{u}]_n)$, take $C$ a compact of $\ddomain$, and define $\nu_C$ such that for $A \in \borel(\ddomain)$
\begin{equation*}
\nu_C (A) := \lambda \int_{g^{-1}(A)} \underset{\delta \in C}{\inf} \frac{ \den (t_\delta(\bs{u})) \dist( \Ldstar[{[\bs{u}]_1}] )^{\lambda-1}}{\dist(\Ld)^\lambda} d\bs{u} \enspace .
\end{equation*}
As the density $\den$ is supposed to be strictly positive on $\R^n$, for all $\delta \in \ddomain$ we have $\dist(\Ld) \geq \dist(\Ld[0]) > 0$.
Using the fact that $\dist$ is a finite measure, and is absolutely continuous, applying the dominated convergence theorem shows that the functions $\delta \mapsto \dist(\Ld)$ and $\delta \mapsto \dist( (-\infty, \delta - [\bs{u}]_1) \times \R^{n-1} \cap \Ld )$ are continuous. Therefore the function $\delta  \mapsto \den(t_\delta (\bs{u})) \dist( \Ldstar[{[\bs{u}]_1}] )^{\lambda-1} / \dist(\Ld)^\lambda $ is continuous and $C$ being a compact, the infimum of this function is reached on $C$ is reached on $C$. Since this function is strictly positive, if $g^{-1}(A)$ has strictly positive Lebesgue measure then $\nu_C(A) > 0$ which proves that this measure is not trivial. By construction $P(\delta,A) \geq \nu_C(A)$ for all $\delta \in C$, so $C$ is a small set which shows that compact sets are small. Since if $\mu_+(A) >0$ we have $P(\delta,A) \geq \nu_C(A) > 0$, the Markov chain $\dmarkov$ is $\mu_+$-irreducible. Finally, if we take $C$ a compact set of $\ddomain$ with strictly positive Lebesgue measure, then it is a small set and $\nu_C(C) > 0$ which means the Markov chain $\dmarkov$ is strongly aperiodic.

The function $\Delta V$ is defined as $\delta mapsto \E(V(\dtt) | \dt = \delta) - V(\delta)$. We want to show a drift condition (see \cite{markovtheory}) on $V$. Using Eq.~\eqref{eq:deltaupdate}
\begin{align*}
\Delta V(\delta) &= \E(\delta - g(\Mtstar) | \dt = \delta) - \delta) \\ 
 &= -\E(g(\Mtstar)) \enspace .
\end{align*}
Therefore using the condition~\eqref{eq:asineq}, we have that there exists a $\epsilon >0$ and a $M \in \R_+$ such that $\forall \delta \in (M,+\infty)$, $\Delta V (\delta ) \leq -\epsilon$. With condtion~\eqref{eq:asint} implies that the function $\Delta V + \epsilon$ is bounded on the compact $[0,M]$ by a constant $b \in \R$. Hence for all $\delta \in \ddomain$
\begin{equation} \label{eq:ineqV}
\frac{\Delta V(\delta)}{\epsilon} \leq -1 + \frac{b}{\epsilon} \ind_{[0,M]}(\delta) \enspace .
\end{equation}
For all $x\in \R$ the level set $C_{V,x}$ of the function $V$, $\lbrace y \in \ddomain | V(y) \leq x \rbrace$, is equal to $[0,x]$ which is a compact set, hence a small set according to what we proved earlier (and hence petite \cite[Proposition 5.5.3]{markovtheory}). Therefore $V$ is unbounded off small sets and with \eqref{eq:ineqV} and Theorem 9.1.8 of \cite{markovtheory}, the Markov chain $\dmarkov$ is Harris recurrent. The set $[0,M]$ is compact and therefore small and petite, so with \eqref{eq:ineqV}, if we denote $r_1$ the constant function $\delta \in \ddomain \mapsto 1$ then with Theorem 14.0.1 of \cite{markovtheory} the Markov chain $\dmarkov$ is positive and is $r_1$-ergodic.

We now want to show a drift condition (see \cite{markovtheory}) on $\Va$.
\begin{align*}
\Delta \Va (\delta) &= \E\left( \exp \left(\alpha \delta - \alpha g\left(\Mtstar\right)\right) | \dt = \delta \right) - \exp\left(\alpha\delta\right) \\
\frac{\Delta \Va}{\Va}(\delta) &= \E\left(\exp\left(-\alpha g\left(\Mtstar\right)\right) | \dt = \delta \right) - 1 \\
 &= \int_{\R^n}  \underset{t \rightarrow +\infty}{\lim}\sum_{k =0}^t \frac{(- \alpha g(\x))^k}{k!} \denstar (\x) d\x -1 \enspace .
\end{align*}
With Eq.~\eqref{eq:denselect} we see that $\denstar(\x) \leq \lambda \den(\x) / \dist(\Ld[0])^\lambda$, so with our assumption that $\E(\exp \alpha | g(\Mtij) | | \dt = \delta ) < \infty $ for $\alpha>0$ small enough we have that the function $\delta \mapsto \E(\exp(\alpha |g(\Mtstar)| | \dt = \delta)$ is bounded for the same $\alpha$. As $\sum_{k=0}^t (- \alpha g(\x))^k/k!\denstar(\x) \leq \exp(\alpha |g(\x)| )\denstar(\x)$ which, with condition \eqref{eq:asexpint}, is integrable so we may apply the theorem of dominated convergence to invert limit and integral:
\begin{align*}
\frac {\Delta \Va}{\Va} (\delta) &=  \underset{t \rightarrow +\infty}{\lim}\sum_{k =0}^t \int_{\R^n} \frac{ (- \alpha g(\x))^k}{k!} \denstar (\x) d\x -1 \\
 &= \sum_{k \in \N} (-\alpha)^k \frac{\E \left({g\left(\Mtstar\right)}^k | \dt = \delta \right)}{k!} - 1 
\end{align*}
Since $\denstar(\x) \leq \lambda \den(\x) / \dist(\Ld[0])^2$, $(-\alpha)^k \E ({g(\Mtstar)}^k | \dt = \delta)/k! \leq (-\alpha)^k \E ({g(\Mtij)}^k)/k!$ which is integrable with respect to the counting measure so we may apply the dominated convergence theorem with the counting measure to invert limit and serie.
\begin{align*}
\underset{\delta \rightarrow +\infty}{\lim} \frac {\Delta \Va}{\Va} (\delta) &=  \sum_{k \in \N} \underset{\delta \rightarrow +\infty}{\lim} \left(-\alpha\right)^k \frac{\E \left({g\left(\Mtstar\right)}^k | \dt = \delta\right)}{k!} - 1 \\
 &= -\alpha \underset{\delta \rightarrow +\infty}{\lim} \E \left(g\left(\Mtstar\right) | \dt=\delta \right) + o\left(\alpha\right) \enspace .
\end{align*}
With condition \eqref{eq:ineqV} we supposed that ${\lim}_{\delta \rightarrow +\infty} \E (g(\Mtstar) | \dt=\delta) > 0$ this implies that for $\alpha>0$ and small enough, $lim_{\delta \rightarrow +\infty} \Delta \Va (\delta)/\Va(\delta) <0$, hence there exists $M \in \ddomain$ and $epsilon >0$ such that $\forall \delta > M$, $\Delta \Va(\delta) < -\epsilon \Va(\delta)$. Finally as $\Delta \Va - \Va$ is bounded on $[0,M]$ there exists $b \in \R$ such that
\begin{equation*}
\Delta \Va (\delta) \leq -\epsilon \Va(\delta) + b \ind_{[0,M]} (\delta) \enspace .
\end{equation*}
According to what we did before in this proof, the compact set $[0,M]$ is small, and hence is petite (\cite[Proposition 5.5.3]{markovtheory}). So the $\mplus$-irreducible Markov chain $\dmarkov$ satisfies the conditions of Theorem 15.0.1 of \cite{markovtheory} which with Theorem 14.0.1 of \cite{markovtheory} proves  that the Markov chain $\dmarkov$ is $\Va$-geometrically ergodic.
\end{myproof}

We now use a law of large numbers (\cite{markovtheory} Theorem~17.0.1) on the Markov chain $(\dt, (\Mtij)_{i\in [1..\lambda], j\in \N})_{t\in \N}$ to obtain an almost sure divergence of the algorithm.

\begin{proposition} \label{pr:lln}
Let a $(1,\lambda)$-ES optimize problem \eqref{eq:pbdef}, handling the constraint through resampling. Assume that the distribution $\dist$ of the random step $\Mtij$ is absolutely continuous with continuous and strictly positive density $\den$, that conditions \eqref{eq:asexpint} and \eqref{eq:asineq} of Proposition~\ref{pr:markovproperties} hold, and denote $\pi$ and $\mu_M$ the stationary distribution of respectively $\dmarkov$ and $\seqmtij$. Then 
\begin{equation} \label{eq:xdiv}
\frac{[\Xt - \Xt[0]]_1}{t}    \overset{a.s.}{\underset{t \rightarrow + \infty}{\longrightarrow}}    \sigma \E_{\pi\times\mu_M}([\Mtstar]_1) \enspace .
\end{equation}

Furthermore if $\E([\Mtstar]_2) < 0$, then the right hand side of Eq.~\eqref{eq:xdiv} is strictly positive.
\end{proposition}

\todo{check $+$ in markov chain updates}
\begin{myproof}
According to Proposition~\ref{pr:markovproperties} the sequence $\dmarkov$ is a Harris recurrent positive Markov chain with invariant measure $\pi$. As $\seqmtij$ is a i.i.d. sequence with distribution $\mu_M$, $(\dt,(\Mtij)_{i \in [1..\lambda], j \in \N})_{t\in \N}$ is also a Harris recurrent positive Markov chain. As $[\Mtstar]_1$ is a function of $\dt$ and $(\Mtij)_{i \in [1..\lambda], j \in \N}$, if $\E_{\pi\times\mu_M}(|[\Mtstar]_1|) < \infty$, according to Theorem 17.0.1 of \cite{markovtheory}, we may apply a law of large numbers on the right hand side of Eq.~\eqref{eq:xrec} to obtain \eqref{eq:xdiv}.

Using Fubini-Tonelli's theorem $\E_{\pi\times\mu_M}(|[\Mtstar]_1|) = \E_{\pi}(\E_{\mu_M}(|[\Mtstar]_1| | \dt = \delta))$. From Eq.~\eqref{eq:denselect} for all $\x \in \R^n$, $\denstar(\x) \leq \lambda \den(\x) / \dist(\Ld[0])^2$, so the condition in \eqref{eq:asexpint} implies that for all $\delta \in \ddomain$, $\E_{\mu_M}(|[\Mtstar]_1| | \dt = \delta)$ is finite. Furthermore, with condition \eqref{eq:asineq}, the function $\delta \in \ddomain \mapsto \E_{\mu_M}(|[\Mtstar]_1| | \dt = \delta)$ is bounded by some $M \in \R$. Therefore as $\pi$ is a probability measure, $\E_\pi(\E_{\mu_M}(|[\Mtstar]_1| | \dt = \delta)) \leq M < \infty$ so we may apply the law of large numbers of Theorem 17.0.1 of \cite{markovtheory}.

Using the fact that $\pi$ is an invariant measure, we have $\E_\pi(\dt) = \E_\pi(\dtt)$, so $\E_\pi(\dt) = \E_\pi(\dt - \sigma g(\Mtstar))$ and hence $\cos \theta \E_\pi([\Mtstar]_1) = - \sin \theta \E_\pi([\Mtstar]_2)$. So using the assumption that $\E([\Mtij]_2) \leq 0$ then we get the strict positivity of $\E_{\pi\times\mu_M}([\Mtij]_1)$.
\end{myproof}

\section{Application to More Specific Distributions}  
\label{sec:examples}        

Throughout this section we give cases where the assumptions on the distribution of the random steps $\dist$ used in Proposition~\ref{pr:markovproperties} or Proposition~\ref{pr:lln} are verified.

The following lemma shows an equivalence between a non-identity covariance matrix for $\dist$ and a different norm and constraint angle $\theta$.

\todo{change dimension index to $k$ in the rest of the paper for consistency}
\begin{lemma} \label{lm:cov}
Let a $(1,\lambda)$-ES optimize problem \eqref{eq:pbdef}, handling the constraint with resampling. Assume that the distribution $\dist$ of the random step $\Mtij$ has positive definite covariance matrix $\C$ with eigenvalues $(\alpha_i^2)_{i \in [1..n]}$ and take $\B = (b_{i,j})_{(i,j) \in [1..n]^2}$ such that $\B \C \B^{-1}$ is diagonal. Denote $\mathcal{A}_{\dist, g, \Xt[0]}$ the sequence of parent points $(\Xt)_{t \in \N}$ of the algorithm with distribution $\dist$ for the random steps $\Mtij$, constraint angle $\theta$ and initial parent $\Xt[0]$.  Then for all $k \in [1..n]$
\begin{equation}
\beta_k \left[\mathcal{A}_{\dist, \theta, \Xt[0]}\right]_k \eqd \left[\mathcal{A}_{\C^{-1/2}\dist, \theta', \Xt[0]'}\right]_k \enspace ,
\end{equation}
where $\beta_k = \sqrt{\sum_{j = 1}^n \frac{b_{j,i}^2}{\alpha_i^2}}$, $\theta' = \arccos (\frac{\beta_1 cos\theta}{\beta_g} )$ with $\beta_g = \sqrt{ \beta_1^2 \cos^2 \theta + \beta_2^2 \sin^2 \theta }$, and $[\Xt[0]']_k = \beta_k[\Xt[0]]_k$ for all $k \in [1..n]$.
\end{lemma}

\begin{myproof}
Take $(\ekb)_{k \in [1..n]}$ the image of $(\ei)_{k \in [1..n]}$ by $\B^{-1}$. We define a new norm $\normb{\cdot}$ such that $\normb{ \ekb } = 1/\alpha_k$. We define two orthonormal basis $(\ekp)_{k \in [1..n]}$ and $(\ekbp)_{k \in [1..n]}$ for $(\R^n, \normb{\cdot})$ by taking $\ekp = \ei[k]/\normb{\ei[k]}$ and $\ekbp = \ekb/\normb{\ekb} = \alpha_k \ekb$.  As $\Var( \Mtij.\ekb ) = \alpha_k^2$, $\Var (\Mtij.\ekbp) = 1$ so in $(\R^n, \normb{\cdot})$ the covariance matrix of $\Mtij$ is the identity.

Take $h$ the function that to $\x \in \R^n$ maps its image in the new orthonormal basis $(\ekp)_{k \in [1..n]}$. As $\ekp = \ei[k]/\normb{\ei[k]}$, $h(\x) = (\normb{\ei[k]} [\x]_k)_{k \in [1..n]}$, where $\normb{\ei[k]} = \normb{ \sum_{i = 1}^n b_{i,k} \ekb } = \sqrt{\sum_{i=1}^n b_{i,k}^2/\alpha_k^2} = \beta_k$. As we changed the norm, the angle between $\nabla f$ and $\nabla g$ is also different in the new space. Indeed $\cos \theta' = h(\nabla g). h(\nabla f)/(\normb{h(\nabla g)}\normb{h(\nabla f)}) = \beta_1^2 \cos \theta / (\sqrt{\beta_1^2 \cos^2 \theta + \beta_2^2 \sin^2 \theta} \beta_1) = \beta_1 \cos \theta / \beta_g$. 

If we take $\Ntij \sim \C^{-1/2} \dist$ then it has the same distribution as $h(\Mtij)$. Take $\Xt' = h(\Xt)$ then for a constraint angle $\theta' = \arccos (\beta_1 \cos \theta / \beta_g)$ and a normalized distance to the constraint $\dt = \Xt'.h(\nabla g) / \st$ the ressampling is the same for $\Ntij$ and $h(\Mtij)$ so $\Nti \eqd h( \Mti)$. Finally the rankings induced by $\nabla f$ or $h(\nabla f)$ are the same so the selection in the same, hence $\Ntstar \eqd h(\Mtstar)$, and therefore $\Xtt' \eqd h(\Xtt)$.
\end{myproof}

Although Eq.~\eqref{eq:xdiv} shows divergence of the algorithm, it is important that it diverges in the right direction, i.e. that the right hand side of Eq.~\eqref{eq:xdiv} has a positive sign. This is achieved when the distribution of the random steps is isotropic, as stated in the following proposition.

\begin{proposition} \label{pr:iso}
Let a $(1,\lambda)$-ES optimize problem \eqref{eq:pbdef} with constant step-size, handling the constraint with resampling. Suppose that the Markov chain $\dmarkov$ is positive Harris, that the distribution $\dist$ of the random step $\Mtij$ is absolutely continuous with strictly positive density $\den$, and take $\C$ its covariance matrix.  If the distribution $\C^{-1/2}\dist$ is isotropic then $\E_{\pi \times \mu_M}([\Mtstar]_1) > 0$.
\end{proposition}

\begin{myproof}
First if $\C = \Id$, using the same method than in the proof of Lemma~\ref{lm:diststep}
\begin{align*}
\denstartwo(y) &= \lambda \int_\R \ldots \int_{\R} \dentilde(u_1, y, u_3, \ldots, u_n) \Pr( u_1 \geq [\Mti]_1)^{\lambda-1} d u_1 \prod_{k=3}^n d u_k \enspace .
\end{align*}
Using Eq.\eqref{eq:denfea} and the fact that the condition $\x \in \Ld$ is equivalent to $[\x]_1 \leq (\delta - [\x]_2 \sin \theta ) / \cos \theta$ we obtain
\begin{equation*}
\denstartwo(y) = \lambda \int_\R \ldots \int_{-\infty}^{\frac{\delta-y \sin\theta}{\cos \theta}} \frac{\den(u_1, y, u_3, \ldots, u_n)}{\dist(\Ld)} \Pr( u_1 \geq [\Mti]_1)^{\lambda-1} d u_1 \prod_{k=3}^n d u_k \enspace .
\end{equation*}
If the distribution of the random steps steps is isotropic then $\den(u_1, y, u_3, \ldots, u_n) = \den(u_1, -y, u_3, \ldots, u_n)$, and as the density $\den$ is supposed strictly positive, for $y >0$ and all $\delta \in $, $\denstartwo(y) - \denstartwo(-y) < 0$ so $\E([\Mtstar]_2 | \dt = \delta) < 0$.
If the Markov chain is Harris recurrent and positive then this imply that $\E_\pi ([\Mtstar]_2) < 0$ and using the reasoning in the proof of Proposition~\ref{pr:lln} $\E_\pi([\Mtstar]_1) > 0$.

For any covariance matrix $\C$ this result is generalized with the use of Lemma~\ref{lm:cov}.
\end{myproof}

Lemma~\ref{lm:cov} and Proposition~\ref{pr:iso} imply the following result to hold for multivariate normal distributions.

\begin{proposition} \label{pr:gauss}
Let a $(1,\lambda)$-ES optimize problem \eqref{eq:pbdef} with constant step-size, handling the constraint with resampling. If $\dist$ is a multivariate normal distribution with mean $\bs{0}$, then $\dmarkov$ is a geometrically ergodic  positive Harris Markov chain, Eq.~\eqref{eq:xdiv} holds and its right hand side is strictly positive.
\end{proposition}

\begin{myproof}
Suppose $\Mtij \sim \Nl(\bs{0}, \Id)$. Then $\dist$ is absolutely continuous and $\den$ is strictly positive. The function $\x \mapsto \exp(g(\x)) \exp(-\|\x\|^2/2)/\sqrt{2\pi}$ is integrable, so Eq.~\eqref{eq:asexpint} is satisfied. Furthermore, when $\delta \rightarrow +\infty$ the constraint disappear so $\Mtij$ behaves like  
$(\Nlambda, \Nl(0,1), \ldots, \Nl(0,1))$ where $\Nlambda$ is the last order statistic of $\lambda$ i.i.d. standard normal variables, so using that $\E(\Nlambda)>0$ and $\E(\Nl(0,1)) = 0$, with multiple uses of the dominated convergence theorem we obtain condition~\eqref{eq:asineq} so with Proposition~\ref{pr:markovproperties} the Markov chain $\dmarkov$ is geometrically ergodic and positive Harris. 

Finally $\dist$ being isotropic the conditions of Proposition~\ref{pr:iso} are fulfilled, and therefore so are every condition of Proposition~\ref{pr:lln} which shows what we wanted.
\end{myproof}

To obtain sufficient conditions for the density of the random steps to be strictly positive, it is advantageous to decompose that distribution into its marginals and the copula combining them. We pay a particular attention to \emph{Archimedean copulas,} i.e., copulas defined 
\begin{gather} 
\label{car}        
                  (\forall \u\in[0,1]^n)\;C_\psi (\u) =\psi
                  (\psi^{-1} ([\u]_1)+\dots+\psi^{-1} ([\u]_n)),  
\end{gather}  
                  where $\psi:[0,+\infty]\to[0,1]$ is an \emph{Archimedean
                  generator}, i.e., $\psi (0)=1,\psi (+\infty )=\lim_{t\to
                  +\infty }\psi (t)=0$, $\psi$ is continuous and strictly
                  decreasing on $[0,\inf\{t:\psi (t)=0 \} )$, and $\psi^{-1}$
                  denotes the generalized inverse of $\psi$,
\begin{gather}
                  (\forall u\in [0,1])\;\psi^{-1} (u)=\inf\{t\in[0,+\infty]:\psi
                  (t)=u \}. 
\end{gather}      
                  The reason for our interest is that Archimedean copulas are
                  invariant with respect to permutations of variables, i.e.,
\begin{gather} 
\label{cip}        
                  (\forall \u\in[0,1]^n)\;C_\psi (\Q\u)=C_\psi (\u).
\end{gather}      
                  holds for any permutation matrix $\Q\in \R^{n,n}$. This
                  can be seen as a weak form of isotropy because in the case of
                  isotropy, (\ref{car}) holds for any rotation matrix, and a
                  permutation matrix is a specific rotation matrix. 

\begin{proposition} 
\label{c3}        
                  Let $\dist$ be the distribution of the two first dimensions of the random step $\Mtij$, $\dist_1$ and $\dist_2$ be its marginals, and $C$ be the copula relating $\dist$ to $\dist_1$ and $\dist_2$. Then the following holds:
\begin{enumerate}
\item             Sufficient for $\dist$ to have a continuous strictly
                  positive density is the simultaneous validity of the following
                  two conditions.
\begin{romanitems}
\item             $\dist_1$ and $\dist_2$ have continuous strictly positive densities
                  $\den_1$ and $\den_2$, respectively.
\item             $C$ has a continuous strictly positive density $c$.
\end{romanitems}
                  Moreover, if (i) and (ii) are valid, then
\begin{gather} 
\label{chg}        
                  (\forall \x\in\mathbb{R}^2)\;h (\x)=c (H_1 ([\x]_1),H_2 ([\x]_2)
                  )h_1([\x]_1)h_2 ([\x]_2). 
\end{gather}      
\item             If $C$ is Archimedean with generator $\psi$, then it is
                  sufficient to replace (ii) with 
\begin{enumerate}
\item[(ii')]      $\psi$ is at least 4-monotone,
                  i.e., $\psi$ is continuous on $[0,+\infty]$, $\psi''$ is
                  decreasing and convex on $\mathbb{R}_+$, and 
                  $(\forall t\in\mathbb{R}_+)\; (-1)^k\psi^{ (k) }
                  (t)\ge 0,k=0,1,2$. 
\end{enumerate}
                  In this case, if (i) and (ii') are valid, then
\begin{gather} 
\label{cha}        
                  (\forall \x \in\mathbb{R}^2)\;h (\x)=\frac{\psi'' (\psi^{-1} (H_1
                  ([\x]_1))+\psi^{-1} (H_2 ([\x]_2)) ) }{\psi' (\psi^{-1} (H_1 ([\x]_1))+\psi^{-1}
                  (H_2([\x]_2)) )}h_1([\x]_1)h_2
                  ([\x]_2).
\end{gather}      
\end{enumerate}   
\end{proposition}

\del{\begin{myproof}
                  The continuity and strict positivity of the density of $H$ is
                  a straightforward consequence of the conditions (i) and (ii),
                  respectively (ii'). In addition, the assumption that $\psi$ is
                  at least 4-monotone implies that it is also 2-monotone, which
                  is for the function $C_\psi$ in (\ref{car}) with $n=2$ a
                  necessary and sufficient condition to be indeed a copula
                  \cite{mcneil09multivariate}. To prove (\ref{chg}), the
                  relationships \small 
\begin{gather} 
\label{cpd}        
                  (\forall \x\in \mathbb{R}^2)\;h (\x)=\frac{\partial^2
                  H}{\partial [\x]_1 \partial [\x]_2 } (\x)
                  ,h_1([\x]_1)=\frac{H_1 }{d[\x]_1} ([\x]_1),h_2([\x]_2)=\frac{H_2
                  }{d[\x]_2} ([x]_2), 
\end{gather}      
                  \normalsize are combined with the Sklar's theorem
                  (\cite{sklar59fonctions}, cf. also \cite{nelsen06introduction}) 
\begin{gather} 
\label{csk}        
                  (\forall \x\in \mathbb{R}^2)\;H (\x)=C (H_1 ([\x]_1),H_2 ([\x]_2))
\end{gather}       
                  and with 
\begin{gather} 
\label{ccd}        
                  c (\bs{u})=\frac{\partial^2 C}{\partial [\bs{u}]_1 \partial [\bs{u}]_2 } (\bs{u}).
\end{gather}      
                  For Archimedean copulas, combining (\ref{ccd}) with
                  (\ref{car}) turns (\ref{chg}) into (\ref{cha}).
\end{myproof}
}

\del{
\begin{example}
                  A simple Archimedean generator is given by  \small 
\begin{gather} 
\label{cgg}        
                  (\forall t\in \mathbb{R}_+)\; \psi (t)=\exp (-t^{ \frac{1}{
                  \vartheta } }) \text{, or equivalently, } (\forall y\in (0,1)
                  )\; \psi^{-1} (y)= (-\ln y)^\vartheta,
\end{gather}      
                  \normalsize with $\vartheta \in[1,+\infty)$. From (\ref{cgg})
                  follows that \small 
\begin{gather} 
\label{cgd}        
                  (\forall t\in \mathbb{R}_+)\; \psi' (t)=-\frac{\exp (-t^{
                  \frac{1}{ \vartheta } })}{ \vartheta t^\frac{ \vartheta-1}{
                  \vartheta } }<0,\psi'' (t)=\frac{\exp (-t^{ \frac{1}{
                  \vartheta } })}{ \vartheta^2 t^\frac{2 \vartheta-1}{ \vartheta
                  } } (t^\frac{1}{\vartheta}+\vartheta-1 )>0
\end{gather}      
                  \normalsize and that the sign of $\psi^{ (k) }$ is opposite to
                  the sign of $\psi^{ (k-1) }$ for $k\ge 2$. Thus $\psi $ is
                  completely monotone ($m$-monotone for every $m\ge 2$), and
                  according to \cite{mcneil09multivariate}, the function defined
                  by (\ref{car}) for this $\psi$, 
\begin{gather} 
\label{gu}        
                  (\forall \bs{u}\in[0,1]^2)\;C_\psi (\bs{u})=\exp \left (-\left ( (-\ln
                  [\bs{u}]_1)^\vartheta+(-\ln [\bs{u}]_2)^\vartheta \right )^{ \frac{1}{
                  \vartheta }}\right  ),  
\end{gather}      
                  indeed is a copula, called \emph{Gumbel copula}
                  \cite{nelsen06introduction}. For this copula, the deinsity $h$
                  can be expressed from (\ref{cha}) using (\ref{cgg}) and
                  (\ref{cgd}). The particular choice $\vartheta=1$ yields the
                  \emph{product copula},
\begin{gather}
                  \Pi(u_1\dots u_m)=\exp (- ( (-\ln u_1)+(-\ln u_2) ) )=u_1u_2, 
\end{gather}      
                  which describes the independence of marginals because due to (\ref{csk}),
\begin{gather}
                  (\forall \x\in \mathbb{R}^2)\;H (x)=H_1 ([\x]_1)H_2 ([\x]_2).
\end{gather}      
\end{example}

}

                  \section{Discussion} \label{sec:end}

\normalsize The paper presents a generalization of recent results of the first author \cite{cah2014linearconstraint} concerning linear optimization by a $(1,\lambda )$-ES in the constant step size case. The generalization consists in replacing the assumption of normality of random steps involved in the evolution strategy by substantially more general distributional assumptions. This generalization shows that isotropic distributions solve the linear problem. Also, although the conditions for the ergodicity of the studied Markov chain accept some heavy-tail distributions, an expnentially vanishing tail allow for geometric ergodicity, which imply a faster convergence to its stationary distribution, and faster convergence of Monte Carlo simulations. In our opinion, these conditions increase the insight into the role that different kinds of distributions play in evolutionary computation, and enlarges the spectrum of possibilities for designing evolutionary algorithms with solid theoretical fundamentals. At the same time, applying the decomposition of a multidimensional distribution into its marginals and the copula combining them, the paper attempts to bring a small contribution to the research into applicability of copulas in evolutionary computation, complementing the more common application of copulas to the Estimation of Distribution Algorithms \cite{cuestainfante10bivariate,salinasgutierez09using,wang12copula}.

Needless to say, more realistic than the constant step size case, but also more difficult to investigate, is the varying step size case. The most important results in \cite{cah2014linearconstraint} actually concern that case. A generalization of those results for non-Gaussian distributions of random steps for cumulative step-size adaptation (\cite{cmaes}) is especially difficult as the evolution path is tailored for Gaussian steps, and some careful tweaking would have to be applied. The $\sigma$ self-adaptation evolution strategy (\cite{beyer1995toward}), studied in \cite{arnold2012behaviour} for the same problem, appears easier, and would be our direction for future research.

\subsection*{     Acknowledgment}

                  The research reported in this paper has been supported by grant ANR-2010-COSI-002 (SIMINOLE) of the French National Research Agency, and Czech Science Foundation (GA\v{C}R) grant 13-17187S.

\bibliographystyle{ieeetr}                   
\bibliography{biblio}

\end{document}